\documentclass[11pt]{amsart}
\usepackage{amsmath,amscd,amssymb,color}
\usepackage[graph,frame,poly,arc,matrix]{xy}  
\RequirePackage{graphicx,fancybox}

\pagespan{1}{6}
\PII{}
\copyrightinfo{}{}

\newcommand{\A}{\mathcal{A}}
\newcommand{\C}{\mathbb{C}}
\newcommand{\R}{\mathbb{R}}
\newcommand{\Q}{\mathbb{Q}}
\newcommand{\Z}{\mathbb{Z}}
\renewcommand{\P}{\mathbb{P}}
\renewcommand{\k}{\Bbbk}

\DeclareMathOperator{\Aut}{Aut}
\DeclareMathOperator{\Hom}{Hom}

\DeclareMathOperator{\gr}{gr}
\DeclareMathOperator{\Tor}{Tor}

\DeclareMathOperator{\Hilb}{Hilb}

\definecolor{brown}{rgb}{0.8,0.1,0.1}
\definecolor{dkblue}{rgb}{0,0.1,0.8}
\definecolor{ltgreen}{rgb}{0.1,0.8,0.1}
\newcommand{\blue}[1]{\textcolor{dkblue}{\bf{#1}}}
\newcommand{\green}[1]{\textcolor{ltgreen}{\bf{#1}}}
\newcommand{\red}[1]{\textcolor{red}{\bf{#1}}}

\begin{document}

\title[Complex hyperplane arrangements]%
{Complex hyperplane arrangements}

\author[M.~Falk]{Michael Falk}
\address{Department of Mathematics and Statistics,  
Northern Arizona University, 
Flagstaff, AZ 86011-5717}
\email{michael.falk@nau.edu} 
\urladdr{http://www.cet.nau.edu/\~{}falk}

\author[A.~I.~Suciu]{Alexander I. Suciu$^*$}
\address{Department of Mathematics,
Northeastern University,
Boston, MA 02115-5000}
\email{a.suciu@neu.edu}
\urladdr{http://www.math.neu.edu/\~{}suciu}
\thanks{$^*$Partially supported by NSF grant DMS-0311142}


\maketitle

In the Fall of 2004, we were fortunate to spend the 
semester in residence at MSRI for the program 
``Hyperplane Arrangements and Applications."  
It was an intense, stimulating, productive, enlightening, 
eventful and most enjoyable experience. It was especially 
so for us ``long-timers" in the field because the program 
truly marked a coming-of-age in the evolution of the 
subject from relative obscurity thirty years ago.  
We had an opportunity to introduce the wonders of
arrangements to a group of graduate students during
the two-week MSRI graduate school in Eugene in early 
August, and to an impressive group of post-docs and 
many other unsuspecting mathematicians during the 
program.  We are glad to have this chance to bring 
some of the ideas to a wider audience.  For general  
reference, we suggest the reader consult the books 
and survey articles listed on the summer school web page, 
\texttt{www.math.neu.edu/\~{}suciu/eugene04.html}. 

In its simplest manifestation, an arrangement is merely a
finite collection of lines in the real plane.  The complement 
of the lines consists of a finite number of polygonal regions.
Determining the number of regions turns out to be a purely
combinatorial problem: one can easily find a recursion for
the number of regions, whose solution is given by a formula
involving only the number of lines and the numbers of lines
through each intersection point. This formula generalizes to
collections of hyperplanes in $\R^\ell$, where the
recursive formula is satisfied by an evaluation of the
characteristic polynomial of the (reverse-ordered) poset of
intersections. The study of characteristic polynomials forms 
the backbone of the combinatorial, and much of the algebraic 
theory of arrangements, which were featured in the MSRI 
workshop ``Combinatorial Aspects of Hyperplane Arrangements" 
last November.

From the topological standpoint, a richer situation is 
presented by arrangements of complex hyperplanes, 
that is, finite collections of hyperplanes in $\C^\ell$ 
(or in projective space $\P^\ell$).  In this case, the 
complement is connected, and its topology, as reflected 
in the fundamental group or the cohomology ring for 
instance, is much more interesting. 

The motivation and many of the applications of the
topological theory arose initially from the connection 
with braids. Let $\A_\ell=\{z_i = z_j \}_{1\leq i < j \leq \ell }$ 
be the arrangement of diagonal hyperplanes in $\C^{\ell}$, 
with complement the configuration space $X_\ell$.  In 1962,  
Fox and Neuwirth showed that $\pi_1(X_\ell)=P_{\ell}$, the 
pure braid group on $\ell$ strings, while Neuwirth and Fadell 
showed that $X_{\ell}$ is aspherical. A few years later, 
as part of his approach to Hilbert's thirteenth problem, 
Arnol'd computed the cohomology ring $H^*(X_\ell,\C)$.

For an arbitrary hyperplane arrangement in $\C^\ell$, 
the fundamental group of the complement,  $G=\pi_1(X)$, 
can be computed algorithmically, using the braid monodromy 
associated to a generic projection of a generic slice in $\C^2$, 
see \cite{CS97} and references there.  The end result is a finite 
presentation with generators $x_i$ corresponding to 
meridians around the $n$ hyperplanes, and commutator 
relators of the form $x_i \alpha_j(x_i)^{-1}$, where 
$\alpha_j\in P_{n}$ are the (pure) braid monodromy 
generators, acting on the meridians via the Artin 
representation $P_n \hookrightarrow \Aut(F_n)$.  

The cohomology ring $H^*(X,\Q)$ was computed 
by Brieskorn in the early 1970's.  His proof 
shows that $X$ is a formal space, 
in the sense of Sullivan:  the rational homotopy 
type of $X$ is determined by $H^*(X,\Q)$.  
In particular, all rational Massey products vanish. 
In 1980, Orlik and Solomon gave a simple 
combinatorial description of the $\k$-algebra 
$H^*(X,\k)$, for any field $\k$:  it is the quotient 
$A=E/I$ of the exterior algebra $E$ on classes 
dual to the meridians, modulo a certain ideal $I$ 
determined by the intersection poset, 
see \cite{OS80, OT92}.

For each $a\in A^1\cong \k^n$, the Orlik-Solomon algebra 
can be turned into a cochain complex $(A,a)$, with 
$i$-th term the degree $i$ graded piece of $A$, 
and with differential given by multiplication by 
$a$, cf.~\cite{Yu95}.  The {\em resonance varieties} 
of $A$ were defined in \cite{Fa97} to be the jumping loci 
for the cohomology of this cochain complex:
\begin{equation}
\label{res-var}
R^{i}_{d}(A)=\{ a\in A^1 \mid \dim_{\k} H^i(A,a) \ge d\}.
\end{equation}

The case of a line arrangement in $\P^2$ is already quite 
fascinating. The subarrangements that contribute components 
to $R^1_d(A)$ have very special combinatorial and geometric 
properties. To be eligible, the incidence matrix for the lines
and intersection points must have null-space of dimension 
at least two, with full support. In addition, the subarrangement 
must have a partition into at least three classes such that no 
point $p$ is incident with one line of one class, while all 
other lines incident with $p$ belong to a second class. 
Such partitions are termed ``neighborly." 
The simplest non-trivial example is provided by the 
braid arrangement $\A_3$, see Figure \ref{fig:braid}.  
In this figure the points represent hyperplanes and the 
lines correspond to the points of multiplicity greater than two. 
This is a diagram of the matroid associated with the arrangement.

\begin{figure}
\setlength{\unitlength}{0.6cm}
\begin{picture}(5,3.5)(0.5,0)
\put(3,3){\line(1,-1){3}}
\put(3,3){\line(-1,-1){3}}
\put(1.5,1.5){\line(3,-1){4.5}}
\put(4.5,1.5){\line(-3,-1){4.5}}
\multiput(0,0)(6,0){2}{\blue{\circle*{0.4}}}
\multiput(1.5,1.5)(3,0){2}{\green{\circle*{0.4}}}
\multiput(3,3)(0,-2){2}{\red{\circle*{0.4}}}
\end{picture}
\caption{\textsf{The $\A_3$ matroid, and a neighborly partition}}
\label{fig:braid}
\end{figure}

When $\k$ has characteristic zero, the Vinberg classification
of generalized Cartan matrices implies an even more
exceptional situation, see \cite{LY00}. One can assign 
multiplicities to the lines so that the partition is into 
classes of equal size $d$, with the same number 
of lines from each class containing each ``inter-class"
intersection point. This partition gives rise to a pencil 
of degree $d$ curves which interpolates the completely 
reducible (not necessarily reduced) curves formed by 
the classes in the partition. The pencil that corresponds 
to Figure \ref{fig:braid} consists of the curves 
$ax^2 + by^2 + cz^2 = 0$, with $a+b+c=0$, see 
Figure \ref{fig:pencil}.  The singular fibers are given by 
$a=0$, $b=0$, and $c=0$.  A non-reduced example is
provided by the arrangement of symmetry planes of the cube
with vertices $(\pm1, \pm 1,\pm 1)$, with the coordinate
hyperplanes having multiplicity two. This multi-arrangement
is interpolated by a pencil of quartics.  Such pencils often
yield  (non-linear) fiberings of the complement by punctured
surfaces, showing in particular that the complement is
aspherical. 

\begin{figure}
\setlength{\unitlength}{0.5cm}
\begin{picture}(6,6)(0,0)
\thicklines
\multiput(2,0)(2,0){2}{\blue{\line(0,1){6}}}
\multiput(0,2)(0,2){2}{\red{\line(1,0){6}}}
\put(0,0){\green{\line(1,1){6}}}
\put(6,0){\green{\line(-1,1){6}}}
\put(3,3){\xy
*\ellipse(7){} 
*\ellipse(10,5.8){} 
\endxy}
\end{picture}
\caption{\textsf{A pencil of conics including 
the braid arrangement}}
\label{fig:pencil}
\end{figure}

There is also an apparent connection between the cohomology 
of $(A,a)$ and critical points of certain multi-variate
rational functions. A resonant degree-one element $a$ is
represented (up to a scalar) by a logarithmic deRham
one-form $d \log \Phi$, where $\Phi$ is a product of the
defining linear forms of the hyperplanes, raised to integral
powers. The dimension of $H^i(A,a)$ is related to the number
of components in the critical locus of $\Phi$ of codimension
$i$. In particular we expect $\Phi$ to have nonisolated
critical points when $a$ is ``generically resonant." 
This is known to be the case for certain high-dimensional
arrangements with certain weights \cite{CV03}, and was
established for arrangements of rank three during the Fall
program.  A  precise description of this relationship in
general is a topic of current study. 

In our example of the $\A_3$ arrangement, $d\log \Phi$ 
is resonant precisely when $\Phi(x,y,z) =
(x^2-y^2)^\alpha(y^2-z^2)^\beta(z^2-x^2)^\gamma$, 
with $\alpha+\beta +\gamma=0$. 
The critical set $d\Phi = 0$ is given (projectively) by 
$[x^2-y^2 : y^2-z^2 : z^2-x^2] = [\alpha:\beta:\gamma]$. 
It is not a coincidence that these critical loci are curves 
in the pencil of Figure \ref{fig:pencil}.

Through the connection with generalized hypergeometric
functions, the critical locus of $\Phi$ is of interest in
relation to the Bethe Ansatz in mathematical physics,  
see \cite{Var03}.  This was a major topic of discussion 
in the MSRI workshop ``Topology of Arrangements and 
Applications" last October.  Somewhat serendipitously, 
the same problem for $\k=\R$ is of interest to the 
combinatorialists studying algebraic statistics, who 
were well-represented in Berkeley last fall.

The {\em characteristic varieties} of a space $X$
are the jumping loci for the cohomology of $X$ with 
coefficients in rank $1$ local systems:
\begin{equation}
\label{char-var}
V^{i}_{d}(X)=\{\mathbf{t}\in 
\Hom (\pi_1(X),\C^*) \mid
\dim_{\C} H^{i}(X,\C_{\mathbf{t}}) \ge d\}, 
\end{equation}
where $\C_{\mathbf{t}}$ denotes the abelian group 
$\C$, with $\pi_1(X)$-module structure given by the 
representation $\mathbf{t}\colon \pi_1(X)\to \C^*$.  

Now suppose $X$ is the complement of an arrangement 
of $n$ hyperplanes. By work of Arapura \cite{Ar97}, the 
irreducible components of the characteristic varieties  
of $X$ are algebraic subtori of the character torus 
$\Hom (\pi_1(X),\C^*)\cong (\C^{*})^n$, 
possibly translated by unitary characters. It turns 
out that the tangent cone at the origin to $V^{i}_{d}(X)$ 
coincides with the resonance variety $R^{i}_{d}(A)$, 
see \cite{CS99, Li01, COr00}. 
Consequently, the resonance varieties are unions 
of linear subspaces; moreover, the algebraic 
subtori in the characteristic varieties are determined 
by the intersection lattice.  Nevertheless, there exist 
arrangements for which the characteristic varieties 
have components that do not pass through the origin, 
\cite{Su01}; it is an open question whether such 
components are combinatorially determined.

Counting certain torsion points on the character torus,
according to their depth with respect to the stratification
by the characteristic varieties, yields information about
the homology of finite abelian covers of the complement, 
see \cite{MS02}.   This approach gives a practical algorithm 
for computing the Betti numbers of the Milnor fiber $F$ of a 
central arrangement in $\C^3$. It has also led to examples 
of multi-arrangements with torsion in $H_1(F)$, see \cite{CDS03}. 
There are no known examples of ordinary arrangements 
with this property.

The tangent-cone theorem, and the linearity of
resonance components, both fail over fields of  
characteristic $p>0$. There is evidence that this 
failure is related to the existence of non-vanishing 
Massey products over $\Z_p$, cf.~\cite{Ma04}. 
In addition, there is an empirical connection between 
translated components of characteristic varieties 
over $\C$ and resonance varieties over fields or 
rings of positive characteristic.  The study of 
resonance varieties in prime characteristic, 
started in \cite{MS00}, leads naturally to the 
theory of line complexes and ruled varieties, 
developed in \cite{Fa04}. The counter-example to the linearity 
question, raised in \cite{Su01}, is a singular, irreducible 
cubic threefold in $\P^4$ ruled by lines, in characteristic 
three. The underlying arrangement is the Hessian arrangement 
of $12$ lines determined by the inflection points on 
a general cubic; see \cite{Fa04}.

As noted by Rybnikov \cite{Ryb98}, the fundamental group 
of the complement, $G=\pi_1(X)$, is not necessarily 
determined by the intersection poset.  Even so, the 
ranks $\phi_k(G)$ of the successive quotients of the 
lower central series $\{G_k\}_{k\ge 1}$, where $G_1=G$ 
and $G_{k+1}=[G,G_{k}]$, are combinatorially 
determined.   Indeed, according to Sullivan, the 
formality of $X$ implies that the graded Lie algebra 
$\gr (G) = \bigoplus_{k\ge 1} G_k/G_{k+1}$ is rationally 
isomorphic to the holonomy Lie algebra $\mathfrak{h}_A$  
associated to $A=H^*(X;\Q)$.  Furthermore, the Chen 
Lie algebra $\gr(G/G''),$ associated to the lower central 
series of $G/G''$, is rationally isomorphic to 
$\mathfrak{h}_A/\mathfrak{h}''_A$, and so the 
Chen ranks $\theta_k(G)$ are also combinatorially 
determined, see \cite{PS04}.  

Much effort has been put in computing explicitly the 
LCS and Chen ranks of an arrangement group $G$.  
It turns out that both can be expressed in terms 
of the Betti numbers of the linear strands in certain  
free resolutions (over $A$ or $E$):
\begin{equation}
\label{lcs-ranks}
\prod_{k=1}^{\infty}(1-t^k)^{\phi_k(G)}=
\sum_{i=0}^{\infty}  \dim \Tor^A_i(\Q,\Q)_i\, t^i , 
\end{equation}
\begin{equation}
\label{chen-ranks}
\theta_k(G)  = 
 \dim \Tor^E_{k-1}(A,\Q)_{k}, \quad \text{for $k\ge 2$}. 
\end{equation}

When the arrangement is of fiber-type (equivalently, 
the intersection lattice is supersolvable), $A$ is a 
Koszul algebra. As noted in \cite{PY99, SY97}, 
formula  \eqref{lcs-ranks} and 
Koszul duality yield the classical LCS formula 
of Kohno \cite{Ko85} and Falk-Randell \cite{FR85}: 
\begin{equation}
\label{lcs}
\prod_{k=1}^{\infty}(1-t^k)^{\phi_k(G)}=\Hilb(A,-t).
\end{equation}

In \cite{Su01}, two conjectures were made, 
expressing (under some conditions) the LCS 
and Chen ranks of an arrangement group in 
terms of the dimensions of the components 
of the first resonance variety.  Write 
$R^1_1(A)=L_1\cup \cdots \cup L_q$, with 
$\dim L_i=d_i$. Then:
\begin{equation}
\label{lcs-res}
\prod_{k=2}^{\infty} (1-t^{k})^{\phi_{k}(G)}=
\prod_{i=1}^{q} \frac{1- d_i t}{(1-t)^{d_i}},
\qquad \text{provided $\phi_4(G)=\theta_4(G)$},
\end{equation}
\begin{equation}
\label{chen-res}
\theta_k(G) = (k-1) \sum_{i=1}^{q}  \binom{k+d_i-2}{k}, 
\qquad \text{for $k$ sufficiently large}.
\end{equation}

The inequality $\ge $ from \eqref{chen-res} has been
proven in \cite{SS05}.  The reverse inequality has 
an algebro-geometric interpretation in terms
of the sheaf on $\C^n$ determined by the linearized
Alexander invariant.  Equality in both \eqref{lcs-res}
and \eqref{chen-res} has been verified for two important 
classes of arrangements:  decomposable arrangements 
(essentially, those for which all components of $R_1^1(A)$ 
arise from sub-arrangements of rank two), and graphic 
arrangements (i.e., sub-arrangements of the braid arrangement); 
see \cite{SS02, PS03, SS05, LS05}. 

Many of the results and observations reported on here 
represent joint work (or work in progress) with our friends 
and collaborators:  Dan Cohen, Graham Denham, 
Dani Matei, Stefan Papadima, Hal Schenck, 
Sasha Varchenko, and Sergey Yuzvinsky. 
Our thanks go to them. In addition, we are grateful 
to many other unnamed participants in the MSRI program 
last Fall, for the countless hours spent in helpful 
and stimulating conversations about arrangements.

\bibliographystyle{amsalpha}

\end{document}